\documentclass{amsart}
\usepackage{latexsym}
\usepackage{amsfonts}
\newcommand{\C}{\mathbb{C}}

\newcommand{\Z}{\mathbb{Z}}
\newcommand{\K}{\mathbf{k}}
\newcommand{\To}{\rightarrow}

\theoremstyle{plain}
\newtheorem{Thm}{Theorem}
\newtheorem{Cor}{Corollary}
\newtheorem{lemma}{Lemma}

\newtheorem{Fact}{Fact}

\theoremstyle{remark}

\newtheorem{remark}{Remark}

\begin{document}
\title[Iterated Images]
{Iterated Images and the Plane Jacobian Conjecture}
\author{Ronen Peretz}
\thanks{Ronen Peretz.
{\small Adress: Dept. of Math., Ben Gurion University of the Negev, Beer-Sheva, 84105, Israel}
{\small Email: ronenp@cs.bgu.ac.il}}
\author{Nguyen Van Chau}\thanks{Nguyen Van Chau.
Supported in part by the National Basic Program of Natural Science, Vietnam, ICTP, Trieste, Italy and ICMC, USP, SP, Brazil. 
{\small Address: Institute of Mathematics, P.O. Box 1078, Hanoi, Vietnam}{\small E-mail: nvchau@math.ac.vn}}
\author{Carlos Gutierrez}
\thanks{Carlos Gutierrez. Supported in part by
PRONEX/FINEP/MCT - grant number  76.97.1080.00.
{ \small Address:  ICMC--USP, P.O.Box 668, 13560--970, Sao Carlos, SP}
{\small  E-mail: gutp@icmc.usp.br}}  
\author{L. Andrew Campbell}\thanks{L. Andrew Campbell. {\small Address: 908 Fire Dance Lane, Pal Desert CA 92211, USA}
{\small Email: landrewcampbell@earthlink.net}}

\include{abstract}

\maketitle
\section{Introduction}

The $n$-dimensional Jacobian conjecture for a field $\K$ of characteristic $0$, 
which we denote
as $JC(\K,n)$, asserts that any polynomial map $f: \K^n \To \K^n$, for which the
Jacobian determinant $J(f)$ is a nonzero constant, is 
bijective and the inverse map is also polynomial (\cite{Survey,ArnoBook}). 
The full Jacobian conjecture asserts that $JC(\K,n)$ is true for all  fields $\K$ 
of characteristic $0$, and all integers $n > 0$; so far, no particular case $JC(\K,n)$
has been proved for any single $\K$ and any single $n > 1$. We consider only the complex 
case $\K = \C$ (it is known that for each $n>0$, $JC(\C,n)$ implies 
$JC(\K,n)$ for any $\K$ of characteristic $0$).

Let $f: \C^n \to \C^n$ be a polynomial map with $J(f)$ a nonzero constant. Then $f$ is an
open map (even locally biholomorphic), $f(\C^n)$ is open and simply connected, 
and $\C^n - f(\C^n)$ is either empty or a closed algebraic subset of $\C^n$ of
codimension at least $2$ (\cite{GaloisCase}). In the case $n=2$ we call $f$ a 
{\em Jacobian pair}.
For a Jacobian pair $f$, the coimage $\C^2 - f(\C^2)$ is empty or at most $0$-dimensional,
hence finite in either case.
We prove a general result about open polynomial endomorphisms with finite coimage,
and then use that result
to provide a recharacterization of $JC(\C,2)$ as a conjecture about certain polynomial 
endomorphisms of cofinite subsets of the plane $\C^2$.

Let $X$ be a closed complex algebraic subset of $\C^N$. 
We equip $X$ with the strong topology that it inherits from the Euclidean topology
on $\C^N$. Let $f: X \To X$ be a polynomial map.
Let $f^k$ denote the iterated map $f^k = f \circ f \circ \dots \circ f$ ($k$ composition
factors).The sequence $X, f(X), f(f(X)), \dots, f^k(X), \ldots$ of iterated images of $X$ 
under $f$ is a nested sequence of sets (each set contains the next or is equal to it).

In general, all we can say about the iterated images of $X$, 
is that they are constructible. That is, they belong to the smallest family of subsets of 
$X$ that contains all the zero sets of polynomials and is closed under complementation
and finite unions and intersections (\cite[II,3.18]{Hartshorne}).

But if $f$ is an open map and $f(X)$ consists of all but finitely many points of $X$,
we will show that the nested sequence $f^k(X)$ eventually stabilizes, with all $f^k(X)$
equal to the same subset of $X$ for $k$ large enough. In addition, we show that if 
$\Omega \subseteq X$ is a cofinite subset of $X$ that is invariant under $f$ 
(that is, $f(\Omega) \subseteq \Omega$), then the iterated
images $f^k(\Omega)$ also stabilize. The second result reduces to the first when 
$\Omega = X$, but it is more general, because, if $\Omega$ is a proper cofinite subset
of $X$, then $\Omega$ is usually not affine 
(e.g. $\C^n$, punctured at one or more points, for $n>1$), hence 
cannot necessarily be embedded as a closed algebraic subset of $\C^N$ for any $N$.
Combining both results into a single theorem yields

\begin{Thm}\label{MainThm} Let  $f: X \To X$ be a polynomial map of a closed complex 
algebraic subset of $\C^N$ to itself. 
If $f$ is an open map, $X-f(X)$ is finite, and $\Omega$ is a cofinite subset of $X$ 
with $f(\Omega) \subseteq \Omega$, 
then there is a $K > 0$ such that
$f^K(\Omega) = f^{K+1}(\Omega) = f^{K+2}(\Omega) = \ldots$ (all $f^k(\Omega)$ are equal for $k \ge K$).
The stable image $\bigcap_{k>0}f^k(\Omega)$ is an open cofinite subset of $\Omega$.
\end{Thm}

We will use the term {\em stability of iterated images}
for the property that the iterated images are all eventually equal.
The restriction of $f$ to $S = \bigcap_{k>0}f^k(\Omega)$
is an open, surjective polynomial map of $S$ to itself.
In a slightly more general context, if $V$ is a complex algebraic variety that can be 
embedded in $\C^N$ as a closed algebraic 
subset, if $f$ is an open
morphism (regular map) of a closed algebraic subset $X$ of $V$ to itself, 
and if $X-f(X)$ is finite,
then stability of iterated images holds for $f$ and any cofinite subset 
$\Omega$ of $X$
satisfying $f(\Omega) \subseteq \Omega$. 
This is true because an embedding of $V$ induces an embedding of $X$,
and $f$ can be expressed as a polynomial map in the coordinates of $\C^N$.

A morphism $f: X \to Y$ of complex algebraic varities is flat if, for every $x \in X$,
the local ring of $X$ at $x$ is a flat module ver the local ring of $Y$ at $f(x)$.
If $f:X \to Y$ is a regular map
of nonsingular complex varieties a geometric definition suffices: $f$ is flat if, and only if,
each fiber $f^{-1}(y)$ of $f$ has dimension $\dim X - \dim Y$
(\cite[App. D.2]{ArnoBook}).
A morphism is \'etale if it is smooth of relative dimension $0$, or
equivalently, flat and unramified (\cite[III, 10.3]{Hartshorne}).
For the case of a regular map of nonsingular complex varieties, these technical conditions
defining the notion of \'etale
reduce to the requirement that the map be locally biholomorphic, i.e. that the
Jacobian determinant of the map (in local coordinates) be nonzero at every point.

\begin{Cor}
Let $n > 1$.
Then stability of iterated images holds for $f: \Omega \To \Omega$, if $f$ is an \'etale 
regular endomorphism 
of a cofinite subset $\Omega$ of $\C^n$ and $f(\Omega)$ is a cofinite subset of $\Omega$.
The stable image of $\Omega$ is an open, simply connected, cofinite subset of $\Omega$.
\end{Cor}

\begin{proof}
$\C^n$ is simply connected; so is the complement in $\C^n$ of any closed complex analytic 
subset of codimension $2$ or more (\cite{LocalAnalytic}).
If $f: \Omega \To \Omega$ is a regular map (a morphism of algebraic varieties), it is a 
rational map
defined at all but finitely many points in $\C^n$, so it is polynomial. 
Let $F$ be the unique extension
of $f$ to a polynomial map from $\C^n$ to $\C^n$. Since $f:\Omega \To \Omega$ is \'etale,
the Jacobian determinant $J(F)$ is nonzero at
every point of $\Omega$. $J(F)$ is a polynomial, so it cannot vanish only on a finite set 
of points.
Since $\Omega$ is cofinite, $J(F)$ must vanish nowhere. Thus $F$ is open. 
Further $\C^n-F(\C^n)$
consists of, at most, some of the finitely many points of $\Omega$ omitted by $f$, 
and some of the
images under $F$ of the finitely many points of $\C^n$ not in $\Omega$. 
So Theorem $1$ applies.
Since the stable image is cofinite in $\Omega$, it is cofinite in $\C^n$, hence simply connected.
\end{proof}

\begin{remark}
In the corollary above, it suffices to assume that $f: \Omega \to \Omega$ has finite
fibers, rather than that it is \'etale. For then $F$ has finite fibers and is dominant
(has a Zariski dense image).
By the fundamental openness principle (\cite[(3.10)]{Mumford}), since $\C^n$ is
nonsingular, $F$ is an open map in the topology we are using (the strong topology
inherited from the Euclidean topology of $\C^n$). Note that $F$ will also be flat,
by the fiber dimension criterion mentioned above for maps between nonsingular varieties
and that all regular flat maps of nonsingular complex varities are open 
(combining the fiber
dimension criterion and the openness principle).
\end{remark}

\begin{Cor}
Let $f$ be an open polynomial map of $\C^2$ to itself. Then the iterated images of $f$
stabilize, and so do the iterated images of any cofinite subset $\Omega \subset \C^2$
such that $f(\Omega) \subseteq \Omega$. The stable image is an open, simply connected,
cofinite subset of $\C^2$ (and of $\Omega$ in the second case).
\end{Cor}

\begin{proof}
By Theorem $1$, it suffices to show that $\C^2 - f(\C^2)$ is finite. It is well
known that $E = \C^2 - f(\C^2)$ is empty or a proper closed algebraic subset of $\C^2$
(see section $3$). Suppose $E$ contains an irreducible component $C$ of dimension $1$.
Let $h = 0$ be a polynomial equation for $C$. Since $C \subseteq E$, 
it follows that $h \circ f$
vanishes nowhere on $\C^2$, so it is a nonzero constant. But then $h$ is constant on the
open set $f(\C^2)$, and hence constant. So $h \equiv c \not= 0 \in \C$ 
and does not vanish on $C$.
We conclude that any irreducible components of $E$ are
of dimension $0$, and so $E$ is finite.
\end{proof}

\begin{remark}
Here, instead of assuming $f$ is open, it suffices to assume that for some
cofinite $\Omega \subseteq \C^2$, the restriction $f|_{\Omega}: \Omega \to \C^2$
has finite fibers and is dominant. In particular, that will be the case if
$f(\Omega) \subseteq \Omega$ and the induced map $\Omega \to \Omega$ has finite
fibers and is dominant.
\end{remark}

As a special case of either of the above two corllaries, we have.

\begin{Cor}
The iterated images of a Jacobian pair stabilize. The stable image is an open, simply connected, cofinite subset of
$\C^2$.
\end{Cor}

It is well known that if $f$ is a Jacobian pair, and $f$ is injective, then $f$ is an automorphism (bijective,
with polynomial inverse). In fact, there is a general principle that an injective endomorphism of a complex
algebraic variety is also surjective (\cite{AxOnInjectivity,IntoOnto}), so bijectivity is easy. 
Moreover, $f$ is an automorphism if its restriction $f|_{\Omega}$ to a dense open
subset of $C^2$ is injective (the birational case of the Jacobian conjecture; this dates back to the origins
of the Jacobian conjecture in \cite{Keller}). So $JC(\C,2)$ holds if every Jacobian pair is injective. 

In contrast, 
little is known about surjectivity of Jacobian pairs or general properties of surjective endomorphisms of
complex algebraic varieties. However, applying  Theorem $\ref{MainThm}$ to Jacobian pairs, we show that

\begin{Thm}
$JC(\C,2)$ is equivalent to the following assertion
\begin{itemize}
\item
If $f$ is a Jacobian pair and $\Omega$ is a cofinite subset of $\C^2$,
such that $f(\Omega) = \Omega$, then $f|_{\Omega}$ is injective.
\end{itemize}
\end{Thm}

\begin{remark} Note that the requirement $f(\Omega) = \Omega$ means that $\Omega$ is invariant
under $f$ and $f|_{\Omega}: \Omega \to \Omega$ is surjective. 
In fact, using reasoning similar to that in Corollary $1$
above, Theorem $2$ can be restated as follows.

{\em The two dimensional complex Jacobian conjecture is true if, and only if, 
every surjective
\'etale endomorphism of a cofinite subset of $\C^2$ is injective.}
\end{remark}

\section{A lemma on the dynamics of a map}
We begin with the following elementary observations on the dynamics of a map of a set to itself.
Neither topology nor algebraic structure are involved.

Let $f: X \To X$ be a map of a nonempty set to itself. Let $f^0 = \text{Id}$, and $f^{k+1} = f \circ f^k$ for $k \ge 0$.
Let $E^k$ denote the set $X - f^k(X)$. Then it is easy to show 
that
$$\emptyset \subseteq E = E^1 \subseteq E^2 \subseteq E^3 \subseteq \ldots .\eqno (*)$$
\noindent and $E^{k+1} \subseteq E^k \cup f(E^k)$.
From  this it follows easily that if $E$ is finite, so are all the $E^k, k > 1$. Also, if, 
for some $K$, $E^K = E^{K+j}$ for some $j > 0$, then all the $E^k$ from $K$ on are the same. 
If there is such a $K$, 
we call the sequence (*) stable. We denote the union of all the $E^k$ by $E^\infty$; 
so (*) is stable if, and only if, $E^K = E^\infty$ for some $K > 0$.

\begin{lemma}\label{dyn-lemma}
Suppose that $E$ is a finite set and (*) is not stable. Then there exists a point
$e \in E$, such that all the points $f^k(e)$ for $k > 0$ belong to $E^\infty$ and are distinct.
For any such $e$ there is
a positive integer $M(e)$,  such that the equation 
$f(x) = f^k(e)$ has exactly one solution for $k \ge M(e)$, namely $x = f^{k-1}(e)$.
\end{lemma}

\begin{proof}
For $a \in X$ define the backward orbit of $a$ as follows. $O(a)$ is a directed graph, whose vertices 
are all points $x \in X$ that satisfy $f^i(x) = a$ for some $i \ge 0$, and in which an edge is directed 
from $x$ to $y$ exactly when $f(y) = x$. Then, if $a \in E^k$, it is clear that $O(a)$ is a tree of 
maximum path length $k-1$. Furthermore, if $a \notin E^\infty$, then $O(a)$ may contain directed 
cycles (arising from periodic points of $f$) and, in any case, has no maximum path length. 

Now assume that (*) is not stable. Then $E$ is nonempty.
We claim that there exists an $e \in E$ such that $f^k(e) \in E^\infty$ for all $k \ge 0$. 
Note that if $x \notin E^\infty$ then $f(x) \notin E^\infty$. So if there is no $e \in E$ for which 
$f^k(e) \in E^\infty$ for all $k \ge 0$, then for each $e \in E$ there is a finite first value $k(e)$ for 
which $f^{k(e)} \notin E^\infty$. Since $E$ is finite, we can define $m = \text{max}_{e \in E}k(e)$. 
If $k > m$ and $a \in E^{k+1} - E^k$ then there is a path of length  $k$ starting at $a \in O(a)$ 
and terminating at some $e \in E$. So $f^k(e) = a \in E^\infty$. That contradicts the definition of $m$. 
So if $k > m$, it follows that $E^{k+1}=E^k$. That contradicts the assumed instability. 
We conclude that some $e \in E$ satisfies $f^k(e) \in E^\infty$ for all $k \ge 0$.
For the remainder of this proof, we fix such an $e$.
Note that all the points $f^k(e)$, $k > 0$ are distinct 
(an equality between any two with different values of $k$ would imply that they are not in $E^\infty$).
 
Consider the directed graph, $\Gamma$, whose set of  vertices is the smallest set that contains $e$ and all the 
points $f^k(e)$, $k > 0$, and contains every $x \in X$ for which $f(x)$ is a vertex. 
Since $f(x) \in E^{k+1}$ implies $x \in E^k$, all the vertices of $\Gamma$ are contained in the graphs 
$O(f^k(e))$. The directed edges of $\Gamma$  go from $x$ to $y$ exactly when $f(y) = x$. 
Let $A$ be the 
set of points $f^k(e)$, $k \ge 0$. 
We claim that the vertices of $\Gamma$ are the disjoint union of $A$ 
and a finite set $B$. For suppose that $v$ is a vertex and $v \notin A$. Then $v$ is $f^m(e')$ for some 
$e' \in E$ that is different from $e$, and some $m \ge 0$. Since $f^{m+j}(e') = f^k(e)$ for some 
$j,k > 0$, only finitely many of the points $f^k(e')$, $k \ge 0$ do not belong to $A$. Since there are
only finitely many possible $e'$, $B$ is finite. For all $k > 0$, the equation $f(x) = f^k(e)$ has the
solution $f^{k-1}(e)$, and that is the only solution in $A$. Any other solution is a vertex of $\Gamma$, 
hence in $B$. Since $B$ is finite, this can occur for only finitely many $k > 0$ (note that inverses of 
distinct elements are distinct). Thus there exists an $M(e) \in \mathbb{N}$ such that $f(x) = f^k(e)$ 
has the unique 
solution $f^{k-1}(e)$ if $k \ge M(e)$.
\end{proof}

\section{Points with a limited number of inverse images}

In this section, let $X$ and $Y$ be closed complex algebraic subsets of $\C^N$, 
and $f: X \To Y$ a polynomial map. For any $n \ge 0$, we define a set $A(f,n)$ as follows.
$$ A(f,n) = \{ y \in Y | \#(f^{-1}(y)) \le n \} $$
\noindent That is, $A(f,n)$ is the set of points of $Y$ which 
have $n$ or fewer inverse images under $f$.
Note that $A(f,0) = E = E^1$ in the notation of the previous section.

The following lemma expresses a well known result. For lack of a suitable
refrence on an elementary level, we provide a simple proof.

\begin{lemma}
If $f: X \To Y$ is an open polynomial map, then each set $A(f,n)$ is a 
closed algebraic subset of $\C^N$ (and of $Y$).
\end{lemma}

\begin{proof}
If we fix $n$ and allow the use of the coefficients of $f$ 
and those of finite sets of polynomials defining $X$ and $Y$,
then there is a first order formula with free variables $z_1,\ldots,z_N$
that is true for precisely those $z = (z_1,\dots,z_N)$ that belong to $A(f,n)$.
So each $A(f,n)$ is constructible: it can be represented as a Boolean combination 
of finitely many closed complex algebraic subsets of $\C^N$. Furthermore, $A(f,n)$ is closed in $\C^N$.
To see this, suppose $y \in Y$ does not belong to $A(f,n)$. Then there are $n+1$
distinct points of $X$ that map to $y$. Since $f: X \To Y$ is open, there is a neighborhood
in $Y$ of $y$ consisting of points with at least $n+1$ inverse images. Thus $A(f,n)$ is a closed
subset of $Y$ and hence of  $\C^N$. But it is classic that a subset of $\C^N$ that is constructible
and closed in the Euclidean topology of $\C^N$ is a closed algebraic subset of $\C^N$ (express the
set as a disjoint union of finitely many sets of the form $C-D$, where $C$ is closed, algebraic, and
irreducible, and $D$ is a proper closed algebraic subset of $C$; then observe that each $C$ is the closure
of $C-D$, hence the set is a finite union of the sets $C$).
\end{proof}

The same holds true if we assume only that $f: X \To Y$ is an open regular map (a morphism) 
of closed complex algebraic subsets of varieties that can be embedded in $\C^N$.
For $X$ and $Y$ can then both be embedded in $\C^N$ for a  large enough $N$, 
and $f$ can be expressed as a polynomial map in the coordinates of $\C^N$.

\section{The iterated image theorem}

For many types of mathematical structures, there is a general principle that an endomorphism that is injective
as a map of sets is also necessarily bijective. Obvious examples are finite sets 
and finite dimensional vector spaces. For complex algebraic sets we have 

\begin{Fact} Let $f: X \To X$ be a polynomial map of a closed algebraic 
subset of $C^N$ to itself and suppose that
$f$ is injective. Then $f$ is surjective.
\end{Fact}

More general statements can be made (\cite{AxOnInjectivity,IntoOnto}) but this will suffice.

\begin{proof}[Proof of Theorem $1$]
Assume that $f: X \To X$ is an open polynomial map of
a closed algebraic subset of $\C^N$ to itself, that $X - f(X)$ is finite, and that $\Omega \subseteq X$ is
a cofinite subset of $X$ with $f(\Omega) \subseteq \Omega$. 

For simplicity, we begin by considering the special case $\Omega = X$. 
Suppose that the sequence (*) is not stable. Then, by Lemma $1$, there
exist $e \in E$ and $M = M(e) > 0$, such that the points $f^k(e)$, $k \ge M$  are all distinct,
they all belong to $E^\infty$, 
and each has a single inverse image under $f$. Let $T = \{f^k(e)\}_{k>M}$. Then $T \subseteq A(f,1)$.
Let $\bar{T}$ be the Zariski closure (in $\C^N$) of $T$; that is, the smallest closed algebraic subset
of $\C^N$ containing $T$. Since $A(f,1)$ is closed algebraic by Lemma $2$, 
$\bar{T} \subseteq A(f,1)$. Obviously $f(T) \subseteq T$ and so, by the Zariski continuity of $f$, we
have $f(\bar{T}) \subseteq \bar{T}$. Thus the restriction of $f$ to $\bar{T}$ is a polynomial map of
$\bar{T}$ to itself. And it is injective because $\bar{T} \subseteq A(f,1)$. So it is surjective,
and $f(\bar{T}) = \bar{T}$. From this it follows that each $t \in \bar{T}$ has an inverse image under $f$
that is also in $\bar{T}$.
Starting with $f^M(e)\in \bar{T}$, we can therefore construct an infinite sequence $\ldots, 
a_{-3}, a_{-2}, a_{-1}, f^M(e)$ in the backwards orbit $O(f^M(e))$, with $f(a_i) = a_{i+1}$, for $i < -1$. 
But that contradicts the fact that $f^M(e) \in E^\infty$. 
We conclude that the sequence (*) must be stable.

Now consider the case in which $\Omega$ is a proper subset of $X$, and suppose that the iterated
images of $\Omega$ do not stabilize. Apply Lemma $1$ to
$f|_{\Omega}: \Omega \To \Omega$. Let $E^k_{\Omega}= \Omega - f^k(\Omega)$ and $E^\infty_{\Omega}
= \bigcup_{k>0}E^k_{\Omega}$. We obtain an $e \in E^1_{\Omega}$, such that the points $f^k(e)$ for $k>0$
are all distinct points of $\Omega$ and lie in $E^\infty_{\Omega}$. 
Also by Lemma $1$, for large enough $k$, each $f^k(e)$ has exactly one inverse image under $f$ that
lies in $\Omega$. If such a point has more than one inverse image under $f$ in $X$, it must be
the image of a point not in $\Omega$. Since there are only finitely many of those, we see that
for large enough $k$, say $k \ge M$, each $f^k(e)$ 
has exactly one inverse image under $f$, namely $f^{k-1}(e)$. Take such a point $f^k(e)$. $k \ge M$.
It belongs to $E^\infty_{\Omega}$. Suppose that $f^k(e) \notin E^\infty$ (computed for $f: X \To X$).
Then there is an infinite path in the backward orbit $O(f^k(e))$. All the points in that path are
distinct (otherwise some points in the forward orbit would not be distinct). If we go back far enough
the points cannot lie in $\Omega$, since $f^k(e) \in E^\infty_{\Omega}$, and hence must be one of
the finitely many points not in $\Omega$. Since the path is infinite, this is a contradiction.
Since we now have all $f^k(e) \in E^\infty \bigcap A(f,1)$ for $k \ge M$, we are in the same 
situation as in the case $\Omega = X$. Just take $T = \{f^k(e)\}_{k>M}$ and proceed as before.
\end{proof}

\section{Iterated images of Jacobian pairs}

Let $f$ be a Jacobian pair, $\Omega$ a subset of $\C^2$ with finite complement, and suppose
$f(\Omega) \subseteq \Omega$. In this situation, we have

\begin{Fact}\label{factZ} If $f|_{\Omega}$ is injective then so is $f$ and 
$f$ is an automorphism.
\end{Fact}

For, since $f$ is injective on a dense open subset of $\C^2$, it must have geometric degree $1$, hence it is
birational and therefore (\cite{Keller,ArnoBook}) an automorphism.

Little is known about surjectivity of Jacobian pairs or their 
$n > 2$ dimensional counterparts. Whether Jacobian pairs are necessarily surjective 
is an open question of long standing (\cite{MR51:5600}).
It is also not known whether a surjective Jacobian pair must be an automorphism (but see
\cite{SurjectiveMaps} for the case in which $f: \C^n \To \C^n$ is a polynomial map with integer coefficients
and $f(\Z^n) = \Z^n$).
Theorem $2$ shows that $JC(\C,2)$ is equivalent to 
the assumption that certain surjective endomorphisms of cofinite subsets of the plane 
are injective; thus seemingly turning on its head the
usual principle that injectivity implies surjectivity.

\begin{proof}[Proof of Theorem $2$]
Let $f$ be a Jacobian pair. If $JC(\C,2)$ holds then $f$ is an automorphism, 
and its
restriction $f|_{\Omega}$ to any subset whatsoever of the plane is injective. Conversely, let us now
assume that $f|_{\Omega}$ is injective if $\Omega$ is any  cofinite subset of $\C^2$ with $f(\Omega) = \Omega$.
Theorem $1$ applies to $f$, so let $\Omega$ be the stable image of $f$. By assumption, $f$ is injective on
$\Omega$, hence  by Fact \ref{factZ}, $f$ is an automorphism. 
\end{proof}

\section{An example}
For concreteness, we illustrate the application of Theorem $1$ with a simple example.
It is taken,  with thanks, from \cite[pp. 292--294]{ArnoBook}, where the properties we list below are
stated. Let $f: \C^2 \To \C^2$
be given by
$$ f(x,y) = (x - 2(xy+1)-y(xy+1)^2, -1 -y(xy+1))$$
Then $f$ is a flat map, hence open, and the image of $\C^2$ under $f$ is $\Omega = \C^2 - (0,0)$.

Since $f(\Omega) = \Omega$, we see that $\Omega$ is the stable image of $f: \C^2 \To \C^2$, as
well as the stable image of $f|_{\Omega}: \Omega \To \Omega$. Note that $f|_{\Omega}$ is an
open, surjective endomorphism of $\C^2 - (0,0)$ that is not injective. It is, of course, not \'etale.

\end{document}